\documentclass[11pt]{article}
\usepackage[margin = 1in]{geometry}
\usepackage{amsmath,epsfig,url,wrapfig}
\usepackage{caption}
\usepackage{subcaption}
\usepackage{bm}
\usepackage{comment}
\usepackage{amssymb}
\usepackage{xcolor}
\usepackage{appendix}
\usepackage{natbib}
\setcitestyle{authoryear,open={(},close={)}}

\newcommand{\innerprod}[1]{\langle\, #1 \,\rangle}

\newcommand{\prob}{\mathbb{P}}
\newcommand{\probz}{P_\mathbf{Z}}

\newtheorem{theorem}{Theorem}

\begin{document}
\baselineskip=16pt
\title{A Gaussian process framework for overlap and causal 
effect estimation with high-dimensional covariates}
\author{Debashis Ghosh \thanks{debashis.ghosh@ucdenver.edu}
\qquad
Efr\'en Cruz Cort\'es \thanks{efren.cruzcortes@gmail.com}}
\date{Department of Biostatistics and Informatics, Colorado School of Public Health, Aurora, CO, 80045}
% this ^ is a cheat. Latex is so bad with author formatting
\maketitle

\begin{abstract}
A powerful tool for the analysis of nonrandomized observational studies has been the potential outcomes model.  Utilization of this framework allows analysts to estimate average treatment effects.  This article considers the situation in which high-dimensional covariates are present and revisits the standard assumptions made in causal inference.  We show that by employing a flexible Gaussian process framework, the assumption of strict overlap leads to very restrictive assumptions about the distribution of covariates, results for which can be characterized using classical results from Gaussian random measures as well as reproducing kernel Hilbert space theory.    In addition, we propose a strategy for data-adaptive causal effect estimation that does not rely on the strict overlap assumption.  These findings reveal the stringency that accompanies the use of the treatment positivity assumption in high-dimensional settings.
\end{abstract}

\newpage

%%% --- INTRODUCTION ---
\section{Introduction}
\label{sec:introduction}
The availability of high-dimensional covariates in administrative databases and electronic health records has led to increasing scientific focus on attempting to evaluate and develop methods for causal inference with these data structures.   There has been a concomitant focus in the statistics and econometrics literature towards the use of machine learning-based methods for performing causal inference with high-dimensional data (e.g., \cite{van2010collaborative}; \cite{van2011targeted}; \cite{athey2016recursive}, \cite{athey2017approximate}, \cite{chernozhukov2018double}).

In light of this work, it is worth revisiting the standard assumptions necessary for performing causal inference in the potential outcomes framework of \cite{rubin1974estimating} and \cite{holland1986statistics}. A key assumption that is needed for proper definition of a causal estimand is the unconfoundedness assumption, which states that the treatment is independent of the potential outcomes conditional on confounders.    Part of the interest in observational studies with high-dimensional covariates is the belief that sufficiently rich sets of covariates will render the unconfoundedness assumption more plausible.    

Another assumption, which is the focus of the current paper, is the treatment positivity assumption, which states that the probability of treatment given covariates is strictly between zero and one for values of the covariate vectors.    This is related to the notion of covariate/confounder overlap between the treatment groups.  This is taken virtually as a given in most causal analyses, but the emergence of clinical decision support systems, deterministic treatment rules and high-dimensional covariates raises the possibility of this assumption being violated.   A simple example of such a violation occurring is a situation in which treatment assignment is made deterministically in a medical setting based on the patient presenting with certain risk factors and comorbidities.  In such a case, the treatment positivity assumption would be violated.  More generally, \cite{robins1997toward} have shown that for estimators of the average causal effect to potentially be semiparametrically efficient, the treatment positivity assumption has to be strengthened to the propensity score being uniformly bounded away from zero and one.   

Methods for performing causal inference with violations of treatment positivity have been limited in the literature. \cite{crump2009dealing} characterized its effects in a setting with limited numbers of covariates and developed a simple rule to exclude subjects based on the propensity score.   The remaining subjects would be those for whom there is sufficient covariate overlap by which one could make valid causal inference.    Traskin and Small (2011) use classification and regression trees (CART) to model group labels derived from the \cite{crump2009dealing} definition of a study population with sufficient overlap to identify factors by which one can define a study population for which one can make causal inferences about.  A recent proposal from \cite{ghosh2017relaxed} suggests defining study populations based on the margin from machine learning algorithms.  Another practical approach is to delete observations with extreme propensities close to zero or one, which is known as propensity score trimming.   Note that these approaches all estimate a causal parameter that is effectively data-dependent; see \cite{ghosh2017relaxed} for further discussion.      

In recent work, \cite{damour2017overlap} studied the notion of strict overlap in the high-dimensional case and described a concept termed strict covariate overlap.   This represents one way of extending the treatment positivity assumption into the high-dimensional setting.    For this situation, they showed that their covariate overlap assumption implied a bound on the discrepancy between the joint distributions of confounders among the treatment and control groups.    The assumption thus places an immediate bound on the difference in joint distributions between treatment and control population.   In order to avoid such a restrictive assumption, \cite{damour2017overlap} suggested that one could make sparsity assumptions of several types.   First, one could assume sparsity at the level of the propensity score.   Second, one could assume the existence of a latent variable that renders treatment independent of the potential outcomes, much as in the standard unconfoundedness assumption.    Finally, one could assume the existence of a low-dimensional subspace such that confounders are independent of the potential outcomes conditional on subspace.   

Central to the notion of strict overlap is that of bounded likelihood ratio, results for which have been developed by \cite{rukhin1993lower,rukhin1997information} and exploited by \cite{damour2017overlap}.    In the current paper, we study the overlap phenomenon using Gaussian processes, which is related to but different from what was done in \cite{damour2017overlap}.   In particular, we will use the theory of random probability measures, and in particular, random Gaussian measures \citep{neveu1968processus,jansson1997gaussian} to characterize implications of overlap.    The line of research we use enjoys a long history in statistical theory, dating back to the results on random measures based on Wiener processes developed by \cite{cameron1944transformations,cameron1945transformations}.    We show that by assuming a Gaussian process framework, the bounded likelihood ratio result presented in \cite{damour2017overlap} leads to an asymptotic implication of overlap that can be characterized by equivalence and orthogonality of Gaussian measures.    These results have been applied to a variety of problems in statistics, including spatial statistics (\cite{stein2012interpolation}) and more recently, classification with functional data (\cite{delaigle2012achieving}, \cite{berrendero2018use}).   The applications of these results to the causal inference 
setting are new.

Next, we show that for a specific Gaussian process model, a phase transition phenomenon exists with respect to the overlap assumption that can be characterized in terms of the component eigenvalues and eigenfunctions.  This will entail the use of and summarizing results from \cite{delaigle2012achieving} and \cite{berrendero2018use}.   The practical implication of these results are manifold:
\begin{enumerate}
    \item For causal inference with high-dimensional covariates, the Gaussian process framework reveals that assuming treatment positivity is equivalent to assuming equivalence of the probability laws of the confounders given treatment groups;
    \item If equivalence does not hold, then one has that the probability laws of the confounders given treatment groups are orthogonal measures, which represents a complete violation of overlap;
    \item We show how one can develop a causal effect estimation strategy that is less reliant on the covariate overlap assumption.  It extends the work of \cite{ghosh2017relaxed} for estimating data-adaptive causal estimands based on the margin.
    
\end{enumerate}
The structure of this paper is as follows.  In Section \ref{sec:background}, we review the potential outcomes model, previous work on covariate overlap as well as provide a review on Gaussian processes.  The latter will prove useful in the framework that we describe in Section \ref{sec:framework}.   In Section \ref{sec:practical}, some practical implications for causal modelling approaches are provided, along with an extension of the strategy from \cite{ghosh2017relaxed} to incorporate covariance structures.  Some discussion concludes the paper in the last section.  

%%% --- BACKGROUND ---
\section{Background}
\label{sec:background}
%%% Preliminaries
\subsection{Preliminaries and Causal inference assumptions}
\label{prelims}
In this paper, we will employ the potential outcomes framework of \cite{rubin1974estimating,holland1986statistics}, which has been widely used in causal modelling. Throughout, we denote the set $\{1, \dots, n\}$ by $[n]$.
We start by assuming an underlying probability space
$\left( \Omega, \Sigma, \prob \right)$, and random variables
$Z: \left( \Omega, \Sigma, \prob \right) \rightarrow \left(\mathcal{Z}, \Sigma_{\mathcal{Z}} \right)$,
$Y: \left( \Omega, \Sigma, \prob \right) \rightarrow \left(\mathcal{Y}, \Sigma_{\mathcal{Y}} \right)$, and
$T: \left( \Omega, \Sigma, \prob \right) \rightarrow \left(\mathbb{B}, \Sigma_{\mathbb{B}} \right)$, where $\mathbb{B}:= \{0,1\}$, and $\mathcal{Z}$, $\mathcal{Y}$ are the confounder and response spaces, respectively.
The derived distribution of $Z$ is defined as $\probz := \prob \circ Z^{-1}$. We define $P_\mathbf{Y}$ and $P_\mathbf{T}$ analogously.
$Y$ denotes the response of interest and ${Z}$ is a $p$-dimensional vector of confounders. $T$ is a binary indicator of treatment exposure, where $T=1$ if treated and $T=0$ if control. We denote their joint distribution by $P$, and assume the observed data, represented as $\{(Y_i,{Z}_i,T_i)\}_{i\in [n]}$, is drawn from $P$. Note that, in this section, we will assume the confounders are simply a finite-dimensional vector.   

Let $\{Y_i(0),Y_i(1)\}$ be the potential outcomes for subject $i$, $i\in [n]$, where $Y_i(0)$ refers to the potential outcome under control and $Y_i(1)$ to that under treatment. We further assume our observation $Y_i$ has the form $Y_i : = Y_i(T_i) = Y_i(1)T_i+ Y_i(0)(1-T_i)$, which is commonly referred to as the consistency assumption in the causal inference literature.   

In an observational study, the vector of covariates ${Z}$ could be related to both the outcome and the treatment assignment. Since both $T$ and the potential outcomes \{Y(0),Y(1)\} are affected by ${Z}$, $ T\perp\{Y(0),Y(1)\}$ will not hold. To enable causal inference in this scenario, we make the following further assumptions.   
\begin{enumerate}
\item Strongly Ignorable Treatment Assumption (SITA):  $\{Y(1),Y(0)\}$ is independent of $T$ given ${Z}$.  
\item Stable Unit Treatment Value Assumption (SUTVA): The potential outcomes for subject $i$ is statistically independent of the potential outcomes for all subjects $j \neq i$, for $i,j \in [n]$.   
\item Treatment Positivity Assumption (TP):  $P_{\mathbf{T} | \mathbf{Z}}(T = 1|{Z}) > 0$ for all values of ${Z}$.
\end{enumerate}
Taking these assumptions in order, SITA means that by conditioning on $Z$, the observed outcomes can be treated as if they come from a randomized complete block design.  \cite{rosenbaum1983central} show that if SITA holds, then the treatment is independent of the potential outcomes given the propensity score, defined as
\begin{align*}
e(Z) := P_{\mathbf{T} | \mathbf{Z}}(T = 1|{Z}).
\end{align*}
\cite{robins1997toward} uses the terminology `no unmeasured confounders' in lieu of SITA.  
The TP assumption was described in the Introduction and will be considered further in the next section.  

In passing, we mention that the typical parameter of focus in causal analyses is the average causal effect (ACE), defined as 
\begin{equation}\label{ace}
ACE :=E(Y(1) - Y(0)),
\end{equation}
where the expectation is taken with respect to the joint distribution of $Y_i(0)$ and $Y_i(1)$.
The use of propensity score modelling \citep{rosenbaum1983central} in conjunction with outcome regression modelling leads to a variety of approaches to average causal effect modelling, a comprehensive overview on the topic being \cite{imbens2015causal}. 

%%% Previous Work
\subsection{Review of previous work regarding treatment positivity and covariate overlap}
\label{sec:prevwork}

We focus on the treatment positivity assumption and describe previous work in this area.   \cite{crump2009dealing} noted the possibility that treatment positivity could be violated and instead defined a subpopulation causal effect using the propensity score.  
Let $\mathbf{I}(C)$ denote the indicator function for the event $C$.   Define the region ${\cal A} := \{ {Z}: c \leq e({Z}) \leq 1 - c \}$ for some $c > 0$. Note we suppress the explicit dependence of ${\cal A}$ on $c$.   \cite{crump2009dealing} define the subpopulation average causal effect as
\begin{align*}
{ACE_{\cal A}} := \frac{\sum_{i=1}^n \mathbf{I}(\{e({Z}_i) \in {\cal A}\}) \left(Y_i(1) - Y_i(0)\right)}{\sum_{i=1}^n \mathbf{I}(\{e({Z}_i) \in {\cal A}\})}.
\end{align*} 
Note the dependence of ${ACE_{\cal A}}$ on the region of the propensity scores that is in ${\cal A}$. In practice, ${\cal A}$ must be estimated from the data, and this is done by estimating the propensity score to obtain $\hat e(Z_i)$, $i=1,\ldots,n$.
By plugging in $\hat e(Z_i)$ instead of $e(Z_i)$ in the definitions of $\mathcal{A}$ and $ACE_{\mathcal{A}}$, we obtain $\widehat{\cal A}$ and $\widehat{ACE}_{{\cal A}}$, respectively. We finally employ $\widehat{ACE}_{\widehat{\cal A}}$ as our final estimator. Based on the variability of the estimated subpopulation average causal effect, \cite{crump2009dealing} proposed an optimization criterion for determining an optimal cutoff value  in the definition of ${\cal A}$ and demonstrated under some mild assumptions that an optimal $c^*$ exists.  The optimal value depends only on the marginal distribution of the propensity scores.  
 
\cite{traskin2011defining} developed a meta-approach for characterizing treatment positivity using classification and regression trees \citep{breiman1984classification}.   In general, a tree classifier works as follows: beginning with a training data set $\{({Z}_i,T_i)\}_{i\in [n]}$, drawn from the joint distribution of $Z$ and $T$, a tree classifier repeatedly splits nodes based on one of the covariates in ${Z}$, until it stops splitting by some stopping criteria (for example, the terminal node only contains training data from one class). Each terminal node is then assigned a class label by the majority of subjects that falls in that terminal node. Once a testing data point with a covariate vector ${Z}$ is introduced, the data point is run from the top of the tree until it reaches one of the terminal nodes. The prediction then will be made by the class label of that terminal node.  Compared to parametric algorithms, tree-based algorithms have several advantages. There is no need to assume any parametric model for a tree; the algorithm for its constructions only requires a criterion for splitting a node and a criterion for when to stop splitting \citep{breiman1984classification}.   \cite{traskin2011defining} propose a modelling approach in which one develops a class label for subject $i$ depending on whether or not $\hat e(Z_i) \in \widehat{\cal A}$, $i\in [n]$. Define $Y^*_i := \mathbf{I}(\hat e(Z_i) \in \widehat{\cal A})$, $i\in [n]$.  \cite{traskin2011defining} then propose fitting a tree model for $Y^*_i$ on $Z_i$ to determine the covariates that explain being in the overlap set of \cite{crump2009dealing}.   We term this a `meta-approach' because the modeling is being based on a response variable that is derived using the estimated propensity score.   

On the theoretical front, \cite{khan2010irregular} demonstrated that if the TP assumption is violated, then irregularities regarding identification and inference about causal effects can occur.  This has echoes in the work of \cite{robins1997toward}, who show that to have regular semiparametric estimators for average causal effects in the high-dimensional case, the model classes for the propensity score and outcome models have to be well-behaved.   Thus, violations in standard overlap assumptions lead to irregularities in estimation and inference.   This was noted in \cite{luo2017estimating}, who found that by assuming a weaker covariate overlap assumption, one could derive an estimator that exhibited super-efficiency (i.e, having an information bound that is smaller than the classical semiparametric information bound for regular estimators).  This phenomenon also occurs in the collaborative targeted maximum likelihood estimator of \cite{van2010collaborative}.   The problem with superefficient estimators, as noted by \cite{damour2017overlap}, is that for model directions where the relaxed assumptions do not hold, the estimators can have unbounded loss functions.   

\cite{damour2017overlap} consider the problem of covariate overlap in one high-dimensional setting. Define the conditional probabilities $P_\ell(Z):= P_{\mathbf{Z} | \mathbf{T}}(Z | T = \ell)$, $\ell\in\mathbb{B}$, and denote their Radon-Nikodym derivatives with respect to $\probz$ as $f_\ell : =dP_\ell / d\probz$. Let $\alpha = P_\mathbf{T}(T = 1)$ and note $\probz = \alpha P_1 + (1-\alpha)P_0$. They note that the strict overlap assumption:
\begin{equation}\label{strictoverlap}
\eta \leq e(Z) \leq 1- \eta \ \ \ \ \ \text{almost everywhere } \probz
\end{equation}
for some $\eta \in (0,0.5)$, is equivalent to the following bounded likelihood ratio assumption:
for $\probz$-almost all $z\in \mathcal{Z}$:
\begin{equation}\label{strictoverlaplr}
\frac{\eta}{1-\eta} \leq \frac{\alpha}{1-\alpha} LR(z) \leq \frac{1- \eta}{\eta}
\end{equation}
where $LR(z): = f_1(z)/f_0(z)$ is the likelihood ratio of the treatment to the control populations. For the sake of completeness, we provide a straightforward proof in the appendix, since none was given in the original paper.
We note assumption (\ref{strictoverlaplr}) is a bounded likelihood ratio assumption. This allows \cite{damour2017overlap} to exploit results from information theory \citep{rukhin1993lower,rukhin1997information} to show that (\ref{strictoverlaplr}) imposes limits on the rate of growth of discriminatory information between the joint distributions of confounders in the treatment and control populations.   Interestingly, these bounds are independent of $p$, the number of confounders.  From an intuitive point of view, this makes sense, for when the number of covariates increase, one would expect that the probability of random classifiers to perfectly separate the data between the treatment and control populations would increase.  This proves problematic for the causal inference problem, as it leads to a violation of the TP assumption.    
%Extending the arguments of D'Amour et al. (2017) to the current setting proves problematic because the likelihood ratio definition in (\ref{strictoverlaplr}) is not well-defined in infinite-dimensional space.   This is because linear and quadratic functions in general are not guaranteed to be well-behaved unless certain topological restrictions are made (e.g., compact linear operators).   

%%% GP Review
\subsection{Review of Gaussian Processes}
\label{sec:gpreview}
We now consider the case where $Z$ is a Gaussian process instead of a finite-dimensional Gaussian variable. $Y$ and $T$ remain as in the previous section but now $Z: \left( \Omega, \Sigma, \prob \right) \rightarrow \left(\mathcal{Z}(\mathcal{T}), \Sigma_{\mathcal{Z}} \right)$ where $\mathcal{Z}(\mathcal{T})$ is a suitable function space on $\mathcal{T}$. To simplify things we can assume these functions are real-valued. It will be more useful to understand $Z$ in the following way: write $Z(\omega,\cdot)$ for the elements of $\mathcal{Z}(\mathcal{T})$, and define the maps $\pi_t (Z):= Z(\cdot,t) \,\,\, \forall t \in \mathcal{T}$. Then, $Z_t : = \pi_t \circ Z$ are functions from $\Omega$ to $\mathbb{R}$. Indeed, each of these functions is a random variable if and only if $Z$ is also a random variable. Indeed $Z$ is fully determined by $\{Z_t\}_{t\in\mathcal{T}}$ and, furthermore, its distribution is fully determined by the finite dimensional joint distributions $\{Z_{t_i}\}_{i\in n}$ for all collections $\{t_i\}_{i\in n}$ and all $n$. When such distributions are all Gaussian, $Z$ is a Gaussian process. If $Z_t$ is square-integrable for all $t\in\mathcal{T}$, we call $Z$ a second-order process, throughout the paper we assume $Z$ is second-order. Comprehensive reviews of Gaussian processes can be found in \cite{neveu1968processus} and \cite{jansson1997gaussian}.
% - RKHS and GPs - 
\subsubsection{Associated Hilbert space}
For simplicity, we assume here that $E(Z_t) = 0$ for all $t\in \mathcal{T}$, although this assumption will be relaxed later in the paper. Gaussian Processes can be characterized by their covariance function, given by 
$$ k(s,t) = \text{Cov}\{Z_s,Z_t\}, \ \ \ \ (s,t) \in {\cal T} \times {\cal T}.$$
For a second-order Gaussian Process, the covariance function $k$ will be a symmetric and positive definite function on ${\cal T} \times {\cal T}$. Well known results from functional analysis ensure each function $k$ of such type is in one-to-one correspondence to a Hilbert space $\mathcal{H}_k$ with the following reproducing property:
$$
\innerprod{f,k(\cdot,t)}_{\mathcal{H}_k} = f(t), \,\,\, \forall f\in\mathcal{H}_t.
$$
In lieu of the reproducing property, We call $\mathcal{H}_k$ the reproducing kernel Hilbert space (RKHS) of $k$, and $k$ the reproducing kernel of $\mathcal{H}_k$. For details see \cite{berlinet2011reproducing, aronszajn1950theory}

%Every Gaussian Process denotes a space of functions on the index set ${\cal T}$, which is known as the Cameron-Martin space.   Let $H$ denote a Gaussian Hilbert space spanned by $\{Z(t)\}$.   For each $\xi \in H$, we can define a function on ${\cal T}$ by $$ R(\xi) \equiv R(\xi)(t) = \innerprod{\xi,Z(t)}_H = E\{\xi Z(t)\}.$$ The Cameron-Martin space is given by $R(H) = \{R(\xi): \xi \in H\}$ and represents a space of real-valued functions on ${\cal T}$. 

%A useful tool for our investigations will be reproducing kernel Hilbert spaces (RKHS) \cite{berlinet2011reproducing}. This is a function space ${\cal H}_K$ that satisfies the property that for any function in it, its pointwise evaluation is a continuous linear functional. As shown in \cite{aronszajn1950theory}, there exists a one-to-one correspondence between ${\cal H}_K$ with a so-called kernel function $K(s,t)$ that is a bounded, symmetric, positive definite function.    The reproducing property of the RKHS states that for any $f \in {\cal H}_K$, $$ \innerprod{f,K(s,\cdot)}_{{\cal H}_K} = f(s).$$   

%The Cameron-Martin space in fact defines an RKHS with reproducing kernel $\rho(s,t)$ and inner product given by $$ \innerprod{f,g}_{R(H)} = \innerprod{R^{-1}f,R^{-1}g}_H.$$ Note that $R^{-1}f$ is well-defined since $R$ is injective. Thus, there is a 1-1 correspondence between the covariance function of a Gaussian process with an RKHS. 

\subsubsection{Equivalence of distributions}
\noindent {\bf Definition.}  Two measures $\mu_0$ and $\mu_1$ defined on a probability space $(X,\mathcal{X})$ are said to be mutually singular (or orthogonal), denoted by $\mu_0 \perp \mu_1$, if there exists a set $A \in {\cal X}$ such that $\mu_0(A) = 0$ and $\mu_1(X \backslash A) = 0$.   $\mu_0$ and $\mu_1$ are said to be mutually equivalent (or equivalent), denoted by $\mu_0 \sim \mu_1$, if for all $A \in {\cal X}$, $\mu_0(A) = 0 \Leftrightarrow \mu_1(A) = 0$.

It is well-known from probability theory that two Gaussian measures defined on the same measure space will either be orthogonal or equivalent \citep{hajek1958property,feldman1958equivalence}.  Assume $\mu$ to be a measure that dominates both $\mu_0$ and $\mu_1$ (e.g., $\mu = \mu_0 + \mu_1$). Define $d\mu_0/d\mu$ and $d\mu_1/d\mu$ to be the Radon-Nikodym derivatives associated with $\mu_0$ and $\mu_1$, respectively.   Furthermore, define 
$$  B := \int \sqrt{\frac{d\mu_1}{d\mu}\frac{d\mu_0}{d\mu}} d\mu$$
and
$$ J := \int \left (\frac{d\mu_0}{d\mu} - \frac{d\mu_1}{d\mu} \right ) \log \left (\frac{d\mu_0/d\mu}{d\mu_1/d\mu} \right ) d\mu.$$
Note that $B$ is the Bhattacharya coefficient between two probability measures and $J$ is the relative entropy.  Several authors have developed characterization results for orthogonality and equivalence of Gaussian measures in terms of $B$ and $J$ \citep{hajek1958property,rao1963discrimination,shepp1966gaussian}   We summarize them here using Theorem 1 from \cite{shepp1966gaussian}.   

\begin{theorem}[Theorem 1 of \cite{shepp1966gaussian}]
\label{th:shepp}
\ \\
(a). $\mu_0 \perp \mu_1$ if and only if $B = 0$ or $J = \infty$.

\noindent (b). $\mu_0 \sim \mu_1$ if and only if $B  > 0$ or $J < \infty$.
\end{theorem}

If $\{\mu_0^{(n)}\}$ and $\{\mu_1^{(n)}\}$ are sequences of Gaussian measures with corresponding Bhattacharya and relative entropy sequences $\{{B}_n\}$ and $\{J_n\}$ converging to $B$ and $J$, then Theorem 1 suggests there are two asymptotic scenarios for Gaussian measures as $n$ approaches infinity. They are either asymptotically orthogonal (part (a)) or equivalent (part (b)). 
%The result is an extension of Kakutani's result for product measures (Kakutani, 1948) as well as an application of zero-one laws in probability (Feller, 1961).  

% Flexibility of Gaussian processes
While the Gaussian process assumption may seem restrictive at first glance, we note that they in fact represent a very flexible class of probability models to fit to data.  They have been widely used in \cite{haran2011gaussian} for spatial statistics, \cite{kennedy2001bayesian} for designs of computer experiments and by \cite{williams2006gaussian} for machine learning.    

%%% --- PROPOSED FRAMEWORK ---
\section{Proposed framework}
\label{sec:framework}
%%% Functional/Causal Assumptions
\subsection{Overlap assumptions and Gaussian Processes}
\label{fcassumptions}
% equivalence versus orthogonal

As before, we represent the random measures for $Z|T = 1$ and $Z|T = 0$ by $P_1$ and $P_0$, respectively. We now reconsider the assumptions of causal inference for $Z$ a Gaussian process.    The functional version of unconfoundedness can be given by 
$$ T \perp \{Y(1),Y(0)\} | {\cal F}_{\mathcal{Z}},$$
where ${\cal F}_{\mathcal{Z}} := \sigma(Z_t: t\in \mathcal{T})$ denotes the smallest $\sigma-$algebra for which every $Z_t$ is measurable.

The assumption of SUTVA remains the same as described in \S 2.  Finally, the treatment positivity assumption can be written as 
$$ 0 < P(T=1|Z) < 1 \text{ almost everywhere} \ \ P_\mathbf{Z}.$$
Analogously to (\ref{strictoverlaplr}) we define strict functional overlap to be 
\begin{equation}\label{fco}
 \frac{\eta}{1-\eta} < \frac{\alpha}{1-\alpha} \frac{dP_1/dP(z)}{dP_0/dP(z)} < \frac{1-\eta}{\eta},
 \end{equation}
for almost-all $z\in\mathcal{Z(\mathcal{T})}$ and for $\eta \in (0, .5)$.   
This is identical to the bounded likelihood ratio assumption of \cite{damour2017overlap}.

Define $\{P^{(n)}_0\}$ and $\{P^{(n)}_1\}$ to be the restrictions of $P_0$ and $P_1$, respectively, to the smallest $\sigma-$algebra for which $\{Z_1,\ldots,Z_n\}$ is measurable.     Let the means be denoted as $\{m^{(n)}_0\}$ and $\{m^{(n)}_1\}$. We have the following result: 

\begin{theorem}
\label{th:sfo_equivalence}
Assumption (\ref{fco}) implies that $P_0 \sim P_1$.  
\end{theorem}

\noindent {\bf Proof:}  Using integration, one can show that   
$$ B(P^{(n)}_0,P^{(n)}_1) = \frac{1}{4} L_n + \frac{1}{8}D^2_n,$$
where 
$$ L_n = 2 \log |\Lambda_n| - \log |\Lambda_{0n}| - \log |\Lambda_{1n}|,$$
with $\Lambda_{jn}$ denoting the covariance operator in $n-$dimensional space for group $j$ $(j=0,1)$ and
$\Lambda_n = (\Lambda_{0n} + \Lambda_{1n})/2$.  The term $D^2_n := \left(m_1^{(n)}- m_0^{(n)}\right)'\Lambda^{-1}_n\left(m_1^{(n)} - m_0^{(n)}\right)$ represents the Mahalanobis distance in $n-$dimensional space.
   By Theorem 3.2 of \cite{rao1963discrimination}, we have that 
\begin{eqnarray*}
 B(P_0,P_1) &=& \lim_{n \rightarrow \infty} B(P^{(n)}_0,P^{(n)}_1)\\
 &=&  \lim_{n \rightarrow \infty} \frac{1}{4} L_n + \frac{1}{8}D^2_n.
 \end{eqnarray*}
As noted in \cite{rao1963discrimination}, $L_n/4$ represents the Hellinger distance between two zero-mean multivariate normal distributions, one with covariance $\Lambda_{0n}$, the other with covariance $\Lambda_{1n}$.  Similarly, $D^2_n/8$ is the Hellinger distance between two $n-dimensional$ multivariate normal distributions with means $m_0^{(n)}$ and $m_1^{(n)}$ and common covariance $\Lambda_n$.   Now, note that $L_n$ and $D_n$ are both greater than equal to zero for all $n$ and are increasing with $n$
The assumption (\ref{fco}), which is a bounded likelihood ratio assumption implies that $L_n$ and $D^2_n$ are uniformly bounded away from zero and infinity for all $n$.   Thus, their limits will also be greater than zero so that $B(P_0,P_1) > 0$.  We can thus use Theorem 1(b) to conclude equivalence.    This concludes the proof of Theorem 2.
\hfill $\blacksquare$

%%% RKHS Interpretation
\subsection{Reinterpretation using RKHS theory and phase transition}

We showed in \S 2.2. how a Gaussian stochastic process can be used to define an RKHS.   In this section, we will assume that $P_0$ and $P_1$ represent the probability laws of Gaussian processes with mean functions $m_0$ and $m_1$, respectively, and shared covariance function $k$. We further assume $m_0 = 0$ and $m_1$, and $k$ are continuous in their respective domains.
%Then, by Theorem 7.1. of \cite{parzen1961approach}, we have that $P_0 \sim P_1$ if and only if $m_1 - m_0 \in {\cal H}_k$, where ${\cal H}_k$ is the RKHS that corresponds to $k$, and the form for the Radon-Nikodym derivative is given by $$ \frac{dP_1}{dP_0}(z) = \exp\left( \innerprod{z,m_1 - m_0}_{k} - 0.5 \innerprod{m_1 - m_0,m_1 - m_0}_{k} \right),$$ where $\innerprod{\cdot,\cdot}_{k}$ denotes the inner product in ${\cal H}_{k}$. Note also that $m_1 - m_0 \notin {\cal H}_k$ if and only if $P_1 \perp P_0$.  As discussed in \S 2.2., with Gaussian processes, the equivalence and orthogonality results represent the only two situations of interest, and for RKHS, it is well-known that the trajectories of $Z$ will be in ${\cal H}_k$ with probability zero or one \citep{lukic2001stochastic}.
If $\mathcal{T}$ is compact, then $k$ is a Mercer kernel and it admits an expansion of the form 
\begin{equation}\label{kl}
 k(s,t) = \sum_{j=1}^{\infty} c_j \psi_j(s) \psi_j(t),   
\end{equation} 
where for $j=1,2,\ldots$, $\{c_j\}$ and $\{\psi_j\}$ denote the eigenvalues and eigenfunctions of the integral operator associated to $k$. For details on such operators, refer to \text{\cite{steinwart2006explicit}}. Based on the decomposition (\ref{kl}), we can decompose the signal $m_1$ as

$$ m_1(t) = \sum_{j=1}^{\infty} a_j \psi_j(t).$$

Recently, \cite{delaigle2012achieving} and \cite{berrendero2018use} have studied and developed this as a model for classification with functional data.   These authors referred to the problem as one of `near-perfect' classification of functional data.  Of relevance is a quote by \cite{delaigle2012achieving}: `these (functional classification) problems have unusual, and fascinating properties that set them apart from their finite-dimensional counterparts.  In particular, we show that in many quite standard settings, the performance of simple (linear) classifiers constructed from training samples becomes perfect as the sizes of the samples diverge....That property never holds for finite-dimensional data, except in pathological cases."  

\cite{delaigle2012achieving} and \cite{berrendero2018use} developed technical characterizations of the `near-perfect' classification phenomenon.   In the context of our discussion, the work of \cite{berrendero2018use} is perhaps the most relevant.  
We have the following theorem.   

\begin{theorem}
\label{th:mercer_restrictions}
\ \\
(a). If (\ref{fco}) holds, then  $\sum_{j=1}^{\infty} c_j^{-1} a_j^2 < \infty$.
%where $$ \psi_j = \sup_{t \in {\cal T}} \psi_j(t).$$

\noindent (b).  If $\sum_{j=1}^{\infty} c_j^{-1} a_j^2 = \infty$, then
(\ref{fco}) cannot hold. 
\end{theorem}

\noindent {\bf Proof:}  We first begin with the proof of Theorem 3a.   By Theorem \ref{th:sfo_equivalence}, if (\ref{fco}) holds, then $P_1 \sim P_0$.  By Theorem 5a of \cite{berrendero2018use}, this is identical to 
$$\sum_{j=1}^{\infty} c_j^{-1} a_j^2 < \infty.$$

\noindent For the proof of Theorem 3b, by Theorem 4b of \cite{berrendero2018use}, $\sum_{j=1}^{\infty} c_j^{-1} a_j^2 = \infty$ is identical to $P_1 \perp P_0$. Note that for probability measures $P_1 \perp P_0$ implies $P_0$ and $P_1$ are not mutually equivalent. Hence, by contraposition of Theorem 2, Assumption (\ref{fco}) cannot hold.
\hfill $\blacksquare$

Theorem~\ref{th:mercer_restrictions} captures a certain type of phase transition in the functional covariate overlap behavior based on the properties of the operator associated to $k$. In particular, the squared eigenvalues of the mean signal divided by the eigenvalues for the kernel operator, characterize strict functional overlap. If strict functional overlap holds, the summation converges (Theorem 3a), while divergence of the series (Theorem 3b) implies strict functional overlap not holding.  
The results of Theorem~\ref{th:mercer_restrictions} are usable in situations where we can derive an explicit Mercer expansion for $k$ in the form of (\ref{kl}). If the  $\{\psi_n\}$ represent a countably infinite set of orthonormal basis functions, then $\{a_n\} \ \ (n \geq 0)$ corresponds to the spectrum of $k$, and if $k$ in (\ref{kl}) is symmetric and positive definite, then $a_n \geq 0$ for all $n$. We now give some examples on the application of Theorem~\ref{th:mercer_restrictions}.

\noindent {\it Example.}  Suppose we set $k$ to be the Gaussian kernel in one dimension: 
$$ k(x,y) = \exp\{-(x - y)^2/2\sigma^2\}, $$
where $x$ and $y$ are scalars, $\sigma^2$ represents a scale parameter.   In \cite{steinwart2006explicit}, the authors show how to represent the Gaussian kernel in the form (\ref{kl}), with $$ c_n = \sqrt{\frac{2\sigma^{2n}}{n}}, \ \ \ \ \psi_n(x) = x^n\exp(-\sigma^2x^2).$$
By Theorem 3a, if (\ref{fco}) holds, then 
$$\sum_{j=1}^{\infty} c_j^{-1} a_j^2 < \infty.$$
In addition, if 
$$\sum_{j=1}^{\infty} c_j^{-1} a_j^2 =  \infty,$$
then (\ref{fco}) cannot hold.  Noting that the summation will converge if $c_j^{-1} a_j^2 = j^{-(1+\epsilon)}$ for some $\epsilon > 0$, some algebra yields that this is equivalent to 
$$ a_j = \sigma^{-j/2}j^{-1/2(1/2+\epsilon)}$$ for $j \geq 0$.  Recalling the definition of $a_j$, this places an immediate constraint on the mean function $m_1(t)$.          
\\
%\noindent {\it Example 2.}  We now let 
%$$ K(x,y) = \sum_{n=0}^\infty a_n \psi_n(x)\psi_n(y),$$
%where $\psi_n(x)$ are the Legendre polynomials, given as solutions to the differential equation
%$$ \frac{d}{dx} \left [(1-x^2)\frac{d\psi_n(x)}{dx} \right ] + n(n+1)\psi_n(x) = 0.$$ 
%Because $\psi = 1$ on the interval $[-1,1]$, where the Legendre polynomials are defined, if $a_n = n$, then by Theorem 3b, (\ref{fco}) cannot hold.  

\section{ Practical Implications for Causal inference}\label{sec:practical}

While the results presented so far have been of a relatively theoretical nature, one practical implication of the results in Section \ref{sec:framework} is that with high-dimensional covariates, there is a strong tension and interplay between (a) maintaining a treatment positivity assumption and (b) having the ability to compute causal effect estimators that are semiparametrically regular in the sense of \cite{robins1997toward}.  One common practice for causal effect estimation is the following strategy:
\begin{enumerate}
    \item Fit a propensity score model for treatment given covariates;
    \item Evaluate covariate balance between the treatment groups, taking into account the propensity scores.  If balance does not hold, modify step 1.
    \item Estimate the causal effect of interest incorporating the propensity score.
\end{enumerate}
If a propensity score model can perfectly predict treatment, then we are in a situation where there the treatment positivity assumption is violated.  Thus, for causal modelling to proceed, one would need stronger assumptions to hold.  For example, targeted learning approaches also fit a model for the potential outcomes, and if that model is correctly specified, then this would overcome the issue of violation of treatment positivity. More fundamentally, the violations in overlap necessitate the use of data-adaptive causal effect estimation approaches, which we describe in the next section.

\subsection{Data-adaptive causal effect estimation}\label{sec:dataadaptive}

   A recent innovation for attempting to relax some of the overlap assumptions is the margin-based approach of \cite{ghosh2017relaxed}.  The intuitive idea of that work is to first identify the margin and then conditional on the margin, develop causal effect estimation and inference procedures.  This is related to, but different from, the finding in \cite{damour2017overlap} that a set of observations where propensity score model misclassifications occur (i.e., the predicted treatment assignment is not concordant with the observed treatment assignment) represent a set of individuals which might be concordant with the treatment positivity assumption.   

We consider the binary treatment scenario and recode treatment as $-1$ and $1$ and use support vector machines (SVM; \cite{cristianini2000introduction}) to fit a propensity score model for treatment.   The objective of SVM is to find a hyperplane, in an appropriate high-dimensional space, that maximizes the margin between the populations defined by $T = 1$ and $T = -1$.   One approach to mapping data into a high-dimensional space is through use of the so-called `kernel trick', which corresponds to using RKHS to define support vector machines.  Using a loss function framework, we can characterize the optimization problem as finding to $f$ to minimize
$$ \sum_{i=1}^n |1 - T_if(Z_i)|_+ \ \ + \lambda\|f\|^2_K,$$
where $|a|_+ = \max(a,0)$, $\lambda > 0$ denotes a smoothing parameter and $\|f\|_K$ denotes the norm of a function $f$ in an RKHS with associated kernel $K$.  At a high level, the approach being proposed here amounts to the following: 
\begin{enumerate}
\item  Fit a support vector machine with kernel $K$ to the data $(T_i,{\bf X}_i)$, $i=1,\ldots,n$; 
\item Determine the observations that are in the margin; denote these as ${\cal M}$; 
\item Estimate the causal effect of interest using $(Y_i,T_i, {\bf X}_i)$, $i \in {\cal M}$;
\end{enumerate}

Extending the arguments of Theorem 1 in \cite{ghosh2017relaxed}, the margin for nonlinear support vector machines based on a kernel $K$ can be viewed as satisfying a functional version of relaxed covariate overlap.  The relaxed covariate overlap notion is detailed more clearly in \cite{ghosh2017relaxed} and as the name suggests, is meant to be a relaxation of the standard covariate overlap assumption in causal inference.

\subsection{Numerical example}\label{sec:rhc}
 
 As an example, we consider the SUPPORT study from \cite{connors1996effectiveness}, a version of which was previously analyzed in \cite{ghosh2015penalized}.   The causal effect of interest is the effect of right heart catherization (RHC) on 30-day survival (dead/alive at 30 days). The dataset contains information on 5735 patients, 2184 of whom received RHC.   Note that the covariates are a mixture of continuous as well as discrete variables, so the assumption of a Gaussian process for the data may not be terribly realistic.  Determining the robustness of results to this assumption is an important topic for future work.  

We employed the implementation of support vector machines as available in the R \textsf{e1071} package based on the Gaussian kernel:
$$ K(x,y) = \exp(-\|x-y\|^2/2\sigma^2),$$
where $x$ and $y$ are $p-$dimensional vectors, and $\sigma >0$ is a scale parameter.
 Out of the original 5735 subjects, 3663 are selected to be in the margin, and based on these subjects, the average causal effect is estimated to be 0.049 with an associated standard error of 0.016, which yields a highly significant effect.  Thus, the use of RHC is associated with decreased probability of 30-day survival, or conversely, increased risk of death within 30 days.  Because of the sensitivity of $K$ to the choice of $\sigma^2$, we performed sensitivity analyses in which we varied $\sigma^2$ and reran the analysis.  In all instances, we found a positive average causal effect.  This aligns with the findings in \cite{connors1996effectiveness}.   

We also used the methodology of \cite{traskin2011defining} to better understand how the margin population was selected.  The plot of the tree is given in Figure 1.    
\begin{figure}
	\centering
	\includegraphics[scale = 0.64]{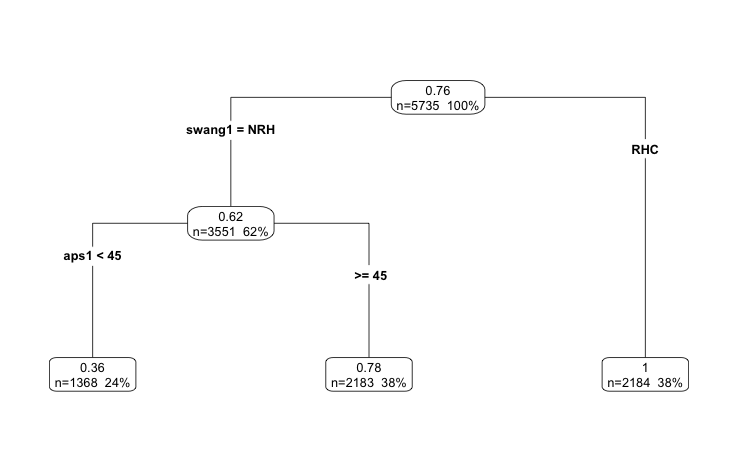}
	\caption{Tree plot from modelling being in the margin as a function of the covariates in the SUPPORT study.  The tree was pruned using a cost-complexity parameter of 0.1.}\label{fig}
\end{figure}
Based on Figure \ref{fig}, we find that all the subjects treated with right-heart catheterization are included in the margin. In addition, patients with higher APACHE scores, which correlate with higher mortality, are much more likely to be included in the margin.  
\section{Discussion}
In this article, we have studied the covariate overlap assumption using Gaussian processes.   The results here complement those in \cite{damour2017overlap}.  The way in which high-dimensional covariates are considered is different from that paper.  There, they consider the situation of $p$, the number of confounders, approaching infinity.  By contrast, here, we use a Gaussian random measure to represent a theoretically infinite-dimensional confounder.   An implication of the analysis here is that the treatment positivity assumption is not as innocuous as it seems, which is similar to the conclusion reached by \cite{damour2017overlap}.  

We also point out that virtually all proposals that are available for causal effect estimation assume a covariate overlap assumption or equivalently, the propensity score being bounded away from zero and one uniformly in the confounders.  Letting $n$ denote the sample size and $p$ the dimension of the confounders, this could be in a `large n, small p' case (e.g., \cite{chan2016globally}) or in a `large p, small  n' setup (e.g., \cite{athey2017approximate}).    Our results suggest that for the latter case, the covariate overlap assumption leads to restrictive assumptions on the distributions of confounders between the treatment groups and that any discrepancy vanishes asymptotically.  If such an assumption is not plausible, then this effectively renders impossible the possibility of average causal effect estimators achieving $n^{1/2}$ convergence \citep{robins1997toward}.   In their words, relaxing the assumption of overlap in the high-dimensional case allows for `pathological' distributions to be considered as part of the model class, and these distributions have the possibility of leading to causal effect estimators with irregular behavior (e.g., superefficiency).   The estimators in \cite{luo2017estimating} and \cite{van2010collaborative} fall into this class.

This article synthesized several results on Gaussian measures that have been in the literature over the last 70 years.   One of the limitations of this article is that the types of covariates that are treated here as being functional are indexed by a one-dimensional set (e.g., time).  An example of biomedical data that could be treated as functional in nature are longitudinal data in electronic medical records, but many covariates in practice might not fit this structure.   Understanding the robustness of the results presented in this paper to violations of the Gaussian process assumption are beyond the scope of the current manuscript and a topic for future research.   

The paper addresses the issue of equivalence versus orthogonality with Gaussian measures for causal inference problems.  In practice, it is difficult to verify if the strict overlap criteria in this article actually hold.   This suggests that approaches that can accommodate violations should be considered.   \cite{petersen2012diagnosing} discussed possible violations of the treatment positivity assumption along with potential causes and tools to evaluate their effects in terms of biases of causal effect estimates.   Suggestions there include restricting the space of treatments, redefining the causal estimand, and using alternative projection functions.   We have proposed redefining the estimand in Section \ref{sec:dataadaptive}. Furthering work in these areas will be necessary in order to advance the use of these analytic approaches for causal effect estimation with high-dimensional confounders.

The recent use of deep learning techniques for learning representations in causal inference has also been proposed 
by \cite{johansson2016learning}.  Extension of the results in Section 3 to the deep learning case would also be of major interest.  \cite{dunlop2018how} have studied `deep' Gaussian processes, which consist of a recursive definition for the Gaussian process based on a Markov chain.   It would be of interest to extend equivalence results to that setting as well.
%If one wishes to model confounders given treatment using a deep Gaussian process, in order to have any hope of doing causal inference, one is effectively assuming that the induced measures are equivalent.  

%To see where it impacts the work of Johansson et al. (2016), we note that they consider discrepancy measures comparing the factual and counterfactual distributions in their optimization problems. Implicitly, to have valid inference, they are assuming that these distributions have common support, which can only occur if the induced measures for confounders given treatment group are identical.   Conversely, if a deep learning algorithm is used to learn the representation and there is, roughly speaking, `separability' between the two distributions, then in an asymptotic limit this leads to the measures being orthogonal.   Thus, one can view this as effectively needing an equivalence in induced measures for valid causal inference to proceed.  

\section*{Acknowledgement}
This research is supported by a pilot grant from the Data Science to Patient Value (D2V) initiative from the University of Colorado.

\appendix
\section*{Appendix}
\subsection*{Proof of Equation (\ref{strictoverlaplr})}
We will prove the right hand side inequality only, as the left hand side inequality is proven analogously. The strict overlap assumption is: $P_{\mathbf{T} | \mathbf{Z}}(T=1 | Z) \leq \eta$ almost everywhere $\probz$. Hence, by Bayes rule:
\begin{align*}
    \frac{P_{\mathbf{Z} | \mathbf{T}}(Z | T = 1) P_{\mathbf{T}}(T = 1)}{\probz(Z)} \leq \eta,
\end{align*}
or, using the notation around equation (\ref{strictoverlaplr}):
\begin{align*}
    \frac{\alpha P_1(Z)}{\alpha P_1(Z) + (1-\alpha)P_0(Z)}\leq \eta \,\,\, & \implies \\
    \alpha P_1(Z) \leq \eta (\alpha P_1(Z) + (1-\alpha) P_0(Z))\,\,\, & \implies \\
    \frac{\alpha P_1(Z)}{(1-\alpha) P_0(Z)} \leq \eta \left(\frac{\alpha P_1(Z)}{(1-\alpha) P_0(Z)} + 1\right) \,\,\, & \implies \\
    \frac{\alpha P_1(Z)}{(1-\alpha) P_0(Z)} (1 - \eta) \leq \eta \,\,\, & \implies \\
    \frac{\alpha P_1(Z)}{(1-\alpha) P_0(Z)} \leq \frac{\eta}{1-\eta}.
\end{align*}
Notice that if a statement holds almost-everywhere $\probz$ then it also holds almost-everywhere $P_0$ and $P_1$ since these are absolutely continuous with respect to $\probz$.

Now, let $b:=\frac{1-\alpha}{\alpha}\frac{\eta}{1-\eta}$. If $P_1(Z)\leq b P_0(Z)$, almost everywhere $\probz$, then $f_1(z) \leq b f_0(z)$ for $\probz$-almost all $z\in \mathcal{Z}$. Otherwise, assuming a set $B$ with positive probability for which $f_1> b f_0$, integrating over $B$ yields $P_1(Z\in B)> b P_0(Z\in B)$, a contradiction.

Conversely, if $f_1(z) \leq b f_0(z)$ for $\probz$-almost all $z\in \mathcal{Z}$, we just need to integrate to obtain $P_1(Z)\leq b P_0(Z)$ almost everywhere $\probz$. Hence, $P_1\leq b P_0$ almost everywhere $\iff f_1 \leq b f_0$ almost everywhere.
\hfill $\blacksquare$

\bibliography{functional_overlap_biblio}

\begin{thebibliography}{43}
\providecommand{\natexlab}[1]{#1}
\providecommand{\url}[1]{\texttt{#1}}
\expandafter\ifx\csname urlstyle\endcsname\relax
  \providecommand{\doi}[1]{doi: #1}\else
  \providecommand{\doi}{doi: \begingroup \urlstyle{rm}\Url}\fi

\bibitem[Aronszajn(1950)]{aronszajn1950theory}
N.~Aronszajn.
\newblock Theory of reproducing kernels.
\newblock \emph{Transactions of the American mathematical society}, 68\penalty0
  (3):\penalty0 337--404, 1950.

\bibitem[Athey and Imbens(2016)]{athey2016recursive}
S.~Athey and G.~Imbens.
\newblock Recursive partitioning for heterogeneous causal effects.
\newblock \emph{Proceedings of the National Academy of Sciences}, 113\penalty0
  (27):\penalty0 7353--7360, 2016.

\bibitem[Athey et~al.(2018)Athey, Imbens, and Wager]{athey2017approximate}
S.~Athey, G.~Imbens, and S.~Wager.
\newblock Approximate residual balancing: De-biased inference of average
  treatment effects in high dimensions.
\newblock \emph{arXiv preprint arXiv:1604.07125}, 2018.

\bibitem[Berlinet and Thomas-Agnan(2011)]{berlinet2011reproducing}
A.~Berlinet and C.~Thomas-Agnan.
\newblock \emph{Reproducing kernel Hilbert spaces in probability and
  statistics}.
\newblock Springer Science \& Business Media, 2011.

\bibitem[Berrendero et~al.(2018)Berrendero, Cuevas, and
  Torrecilla]{berrendero2018use}
J.~R. Berrendero, A.~Cuevas, and J.~L. Torrecilla.
\newblock On the use of reproducing kernel hilbert spaces in functional
  classification.
\newblock \emph{Journal of the American Statistical Association}, 113\penalty0
  (523):\penalty0 1210--1218, 2018.

\bibitem[Breiman et~al.(1984)Breiman, Friedman, Olshen, and
  Stone]{breiman1984classification}
L.~Breiman, J.~Friedman, R.~Olshen, and C.~Stone.
\newblock \emph{Classification and regression trees.}
\newblock Statistics/Probability Series. Wadsworth \& Brooks/Cole, 1984.

\bibitem[Cameron and Martin(1945)]{cameron1945transformations}
R.~Cameron and W.~Martin.
\newblock Transformations of wiener integrals under a general class of linear
  transformations.
\newblock \emph{Transactions of the American Mathematical Society}, 58\penalty0
  (2):\penalty0 184--219, 1945.

\bibitem[Cameron and Martin(1944)]{cameron1944transformations}
R.~H. Cameron and W.~T. Martin.
\newblock Transformations of wiener integrals under translations.
\newblock \emph{Annals of Mathematics}, pages 386--396, 1944.

\bibitem[Chan et~al.(2016)Chan, Yam, and Zhang]{chan2016globally}
K.~C.~G. Chan, S.~C.~P. Yam, and Z.~Zhang.
\newblock Globally efficient non-parametric inference of average treatment
  effects by empirical balancing calibration weighting.
\newblock \emph{Journal of the Royal Statistical Society: Series B (Statistical
  Methodology)}, 78\penalty0 (3):\penalty0 673--700, 2016.

\bibitem[Chernozhukov et~al.(2018)Chernozhukov, Chetverikov, Demirer, Duflo,
  Hansen, Newey, and Robins]{chernozhukov2018double}
V.~Chernozhukov, D.~Chetverikov, M.~Demirer, E.~Duflo, C.~Hansen, W.~Newey, and
  J.~Robins.
\newblock Double/debiased machine learning for treatment and structural
  parameters.
\newblock \emph{The Econometrics Journal}, 2018.
\newblock Accepted Author Manuscript. doi:10.1111/ectj.12097.

\bibitem[Connors et~al.(1996)Connors, Speroff, Dawson, Thomas, Harrell, Wagner,
  Desbiens, Goldman, Wu, Califf, et~al.]{connors1996effectiveness}
A.~F. Connors, T.~Speroff, N.~V. Dawson, C.~Thomas, F.~E. Harrell, D.~Wagner,
  N.~Desbiens, L.~Goldman, A.~W. Wu, R.~M. Califf, et~al.
\newblock The effectiveness of right heart catheterization in the initial care
  of critically iii patients.
\newblock \emph{Jama}, 276\penalty0 (11):\penalty0 889--897, 1996.

\bibitem[Cristianini et~al.(2000)Cristianini, Shawe-Taylor,
  et~al.]{cristianini2000introduction}
N.~Cristianini, J.~Shawe-Taylor, et~al.
\newblock \emph{An introduction to support vector machines and other
  kernel-based learning methods}.
\newblock Cambridge university press, 2000.

\bibitem[Crump et~al.(2009)Crump, Hotz, Imbens, and Mitnik]{crump2009dealing}
R.~K. Crump, V.~J. Hotz, G.~W. Imbens, and O.~A. Mitnik.
\newblock Dealing with limited overlap in estimation of average treatment
  effects.
\newblock \emph{Biometrika}, 96\penalty0 (1):\penalty0 187--199, 2009.

\bibitem[D'Amour et~al.(2017)D'Amour, Ding, Feller, Lei, and
  Sekhon]{damour2017overlap}
A.~D'Amour, P.~Ding, A.~Feller, L.~Lei, and J.~Sekhon.
\newblock Overlap in observational studies with high-dimensional covariates.
\newblock \emph{arXiv preprint arXiv:1711.02582}, 2017.

\bibitem[Delaigle and Hall(2012)]{delaigle2012achieving}
A.~Delaigle and P.~Hall.
\newblock Achieving near perfect classification for functional data.
\newblock \emph{Journal of the Royal Statistical Society: Series B (Statistical
  Methodology)}, 74\penalty0 (2):\penalty0 267--286, 2012.

\bibitem[Dunlop et~al.(2018)Dunlop, Girolami, Stuart, and
  Teckentrup]{dunlop2018how}
M.~M. Dunlop, M.~A. Girolami, A.~M. Stuart, and A.~L. Teckentrup.
\newblock How deep are deep gaussian processes?
\newblock \emph{Journal of Machine Learning Research}, 19\penalty0
  (54):\penalty0 1--46, 2018.
\newblock URL \url{http://jmlr.org/papers/v19/18-015.html}.

\bibitem[Feldman(1958)]{feldman1958equivalence}
J.~Feldman.
\newblock Equivalence and perpendicularity of gaussian processes.
\newblock \emph{Pacific Journal of Mathematics}, 8\penalty0 (5):\penalty0
  699--708, 1958.

\bibitem[Ghosh(2018)]{ghosh2017relaxed}
D.~Ghosh.
\newblock Relaxed covariate overlap and margin-based causal effect estiamtion.
\newblock \emph{Statistics in Medicine}, 37\penalty0 (28):\penalty0 4252--4265,
  2018.
\newblock available at arxiv.org/abs/1801.00816.

\bibitem[Ghosh et~al.(2015)Ghosh, Zhu, and Coffman]{ghosh2015penalized}
D.~Ghosh, Y.~Zhu, and D.~L. Coffman.
\newblock Penalized regression procedures for variable selection in the
  potential outcomes framework.
\newblock \emph{Statistics in Medicine}, 34\penalty0 (10):\penalty0 1645--1658,
  2015.

\bibitem[Hajek(1958)]{hajek1958property}
J.~Hajek.
\newblock A property of j-divergences of marginal probability distributions.
\newblock \emph{Czechoslovak Mathematical Journal}, 8\penalty0 (3):\penalty0
  460--462, 1958.

\bibitem[Haran(2011)]{haran2011gaussian}
M.~Haran.
\newblock Gaussian random field models for spatial data.
\newblock In S.~Brooks, A.~Gelman, G.~Jones, and X.~Meng, editors,
  \emph{Handbook of {M}arkov chain {M}onte {C}arlo}. Springer-Verlag, 2011.

\bibitem[Holland(1986)]{holland1986statistics}
P.~Holland.
\newblock Statistics and causal inference (with discussion).
\newblock \emph{Journal of American Statistical Association}, 81\penalty0
  (396):\penalty0 945--970, 1986.

\bibitem[Imbens and Rubin(2015)]{imbens2015causal}
G.~W. Imbens and D.~B. Rubin.
\newblock \emph{Causal inference in statistics, social, and biomedical
  sciences}.
\newblock Cambridge University Press, 2015.

\bibitem[Jansson(1997)]{jansson1997gaussian}
S.~Jansson.
\newblock \emph{Gaussian Hilbert Spaces}.
\newblock Cambridge Tracts in Mathematics. Cambridge, 1997.

\bibitem[Johansson et~al.(2016)Johansson, Shalit, and
  Sontag]{johansson2016learning}
F.~Johansson, U.~Shalit, and D.~Sontag.
\newblock Learning representations for counterfactual inference.
\newblock \emph{ICML}, pages 3020--3029, 2016.

\bibitem[Kennedy and O'Hagan(2001)]{kennedy2001bayesian}
M.~Kennedy and A.~O'Hagan.
\newblock Bayesian calibration of computer models.
\newblock \emph{Journal of the Royal Statistical Society: Series B (Statistical
  Methodology)}, 63\penalty0 (3):\penalty0 425--464, 2001.

\bibitem[Khan and Tamer(2010)]{khan2010irregular}
S.~Khan and E.~Tamer.
\newblock Irregular identification, support conditions and inverse weight
  estimation.
\newblock \emph{Econometrica}, 78\penalty0 (6):\penalty0 2021--2042, 2010.

\bibitem[Luo et~al.(2017)Luo, Zhu, and Ghosh]{luo2017estimating}
W.~Luo, Y.~Zhu, and D.~Ghosh.
\newblock On estimating regression causal effects using sufficient dimension
  reduction.
\newblock \emph{Biometrika}, 104\penalty0 (1):\penalty0 51--65, 2017.

\bibitem[Neveu(1968)]{neveu1968processus}
J.~Neveu.
\newblock Processus al\'eatoire gaussien.
\newblock \emph{Séminaire de mathématiques supérieures}, 1968.

\bibitem[Petersen et~al.(2012)Petersen, Porter, Gruber, Wang, and van~der
  Laan]{petersen2012diagnosing}
M.~L. Petersen, K.~E. Porter, S.~Gruber, Y.~Wang, and M.~J. van~der Laan.
\newblock Diagnosing and responding to violations in the positivity assumption.
\newblock \emph{Statistical methods in medical research}, 21\penalty0
  (1):\penalty0 31--54, 2012.

\bibitem[Rao and Varadarajan(1963)]{rao1963discrimination}
C.~R. Rao and V.~Varadarajan.
\newblock Discrimination of gaussian processes.
\newblock \emph{Sankhy{\=a}: The Indian Journal of Statistics, Series A}, pages
  303--330, 1963.

\bibitem[Robins and Ritov(1997)]{robins1997toward}
J.~M. Robins and Y.~Ritov.
\newblock Toward a curse of dimensionality appropriate (coda) asymptotic theory
  for semi-parametric models.
\newblock \emph{Statistics in medicine}, 16\penalty0 (3):\penalty0 285--319,
  1997.

\bibitem[Rosenbaum and Rubin(1983)]{rosenbaum1983central}
P.~R. Rosenbaum and D.~B. Rubin.
\newblock The central role of the propensity score in observational studies for
  causal effects.
\newblock \emph{Biometrika}, 70\penalty0 (1):\penalty0 41--55, 1983.

\bibitem[Rubin(1974)]{rubin1974estimating}
D.~B. Rubin.
\newblock Estimating causal effects of treatments in randomized and
  nonrandomized studies.
\newblock \emph{Journal of educational Psychology}, 66\penalty0 (5):\penalty0
  688, 1974.

\bibitem[Rukhin(1997)]{rukhin1997information}
A.~Rukhin.
\newblock Information-type divergence when the likelihood ratios are bounded.
\newblock \emph{Applicationes Mathematicae}, 24:\penalty0 415--423, 1997.

\bibitem[Rukhin(1993)]{rukhin1993lower}
A.~L. Rukhin.
\newblock Lower bound on the error probability for families with bounded
  likelihood ratios.
\newblock \emph{Proceedings of the American Mathematical Society}, 119\penalty0
  (4):\penalty0 1307--1314, 1993.

\bibitem[Shepp(1966)]{shepp1966gaussian}
L.~Shepp.
\newblock Gaussian measures in function space.
\newblock \emph{Pacific Journal of Mathematics}, 17\penalty0 (1):\penalty0
  167--173, 1966.

\bibitem[Stein(2012)]{stein2012interpolation}
M.~L. Stein.
\newblock \emph{Interpolation of spatial data: some theory for kriging}.
\newblock Springer Science \& Business Media, 2012.

\bibitem[Steinwart et~al.(2006)Steinwart, Hush, and
  Scovel]{steinwart2006explicit}
I.~Steinwart, D.~Hush, and C.~Scovel.
\newblock An explicit description of the reproducing kernel hilbert spaces of
  gaussian rbf kernels.
\newblock \emph{IEEE Transactions on Information Theory}, 52\penalty0
  (10):\penalty0 4635--4643, 2006.

\bibitem[Traskin and Small(2011)]{traskin2011defining}
M.~Traskin and D.~S. Small.
\newblock Defining the study population for an observational study to ensure
  sufficient overlap: a tree approach.
\newblock \emph{Statistics in Biosciences}, 3\penalty0 (1):\penalty0 94--118,
  2011.

\bibitem[van~der Laan and Gruber(2010)]{van2010collaborative}
M.~J. van~der Laan and S.~Gruber.
\newblock Collaborative double robust targeted maximum likelihood estimation.
\newblock \emph{The international journal of biostatistics}, 6\penalty0 (1),
  2010.

\bibitem[Van~der Laan and Rose(2011)]{van2011targeted}
M.~J. Van~der Laan and S.~Rose.
\newblock \emph{Targeted learning: causal inference for observational and
  experimental data}.
\newblock Springer Science \& Business Media, 2011.

\bibitem[Williams and Rasmussen(2006)]{williams2006gaussian}
C.~K. Williams and C.~E. Rasmussen.
\newblock \emph{Gaussian processes for machine learning}.
\newblock MIT Press Cambridge, MA, 2006.

\end{thebibliography}

%\printbibliography
\end{document}